\documentclass[9pt,a4paper]{extarticle}

\usepackage{amssymb}
\usepackage{amsmath}
\usepackage{amsthm}
\usepackage[numbers]{natbib}
\usepackage{microtype}
\usepackage[final]{showkeys}
\usepackage{graphicx}
\usepackage{bm,color}

\newtheorem{Theorem}{Theorem}[section]
\newtheorem{Remark}{Remark}[section]
\newtheorem{Definition}{Definition}[section]
\newtheorem{Corollary}{Corollary}[section]
\newtheorem{ex}{Example}[section]

\usepackage[hmargin=3cm,vmargin={3.5cm,4cm}]{geometry}
\numberwithin{equation}{section}
\allowdisplaybreaks

\title{\textbf{On Some Operators Involving Hadamard Derivatives}}

\begin{document}

	\author{$\text{Roberto Garra}_1$, $\text{Federico Polito}_2$\\
		\footnotesize (1) -- Dipartimento di Scienze di Base ed Applicate per l'Ingegneria, \\
		\footnotesize ``Sapienza'' Universit\`a di Roma.\\
		\footnotesize Via A. Scarpa 16, 00161, Rome, Italy.\\
		\footnotesize Email address: roberto.garra@sbai.uniroma1.it \\
		\footnotesize (2) -- Dipartimento di Matematica, Universit\`a di Torino.\\
		\footnotesize Via Carlo Alberto 10, 10123, Torino, Italy.\\
		\footnotesize Tel: +39-011-6702937, fax: +39-011-6702878\\
		\footnotesize Email address: federico.polito@unito.it
		}
	
	\maketitle

	\begin{abstract} 
		\noindent	
		In this paper we introduce a novel Mittag--Leffler-type function
		and study its properties in relation to some integro-differential operators involving
		Hadamard fractional derivatives or Hyper-Bessel-type operators.
		We discuss then the utility of these results
		to solve some integro-differential equations involving these operators
		by means of operational methods.
		We show the advantage
		of our approach through some examples. Among these, an application to a modified Lamb--Bateman integral equation
		is presented.
		
		\vspace{.2cm}
	\end{abstract}
	
	\section{Introduction}

		Fractional calculus is a developing field of the applied analysis concerning methods and
		applications of integro-differential equations
		involving fractional operators. In literature there are  many different definitions
		of fractional derivatives and integrals  \citep{libro}.
		One of these definitions was introduced by \citet{Hadamard} in 1892 and, although not so frequently
		used in applications, it is generally mentioned in classical reference books on fractional calculus
		\citep[Section 18.3]{libro}.
		In some recent papers by \citet{Babusci,Babusci1}, it has been shown by means of
		operational techniques that the Hadamard integral 
		corresponds to a negative fractional power of the operator $x \frac{\mathrm d}{\mathrm dx}$. By
		exploiting the integral representation of such operator, they have shown that
		the Lamb--Bateman integral equation can be rearranged as a fractional integro-differential
		equation involving Hadamard fractional derivatives of order $1/2$.
		However, a very few works on applications of Hadamard integrals
		and derivatives have been carried out, mainly due to the intrinsic intractability of such operators.

		On the other hand, many papers on the properties of Laguerre derivatives and their
		applications in partial differential equations (see for 
		example \citet{Dattoli}) are present in literature. By means of operational methods \citet{Dattoli1}
		have found a simple way to solve analitically
		a class of generalized evolution problems involving Laguerre-type operators.
		The purpose of this paper is to extend these results and to consider some applications in order to
		solve integro-differential equations
		involving Hadamard fractional integrals and derivatives.
		
		First, we define a general class of functions related to fractional operators involving Hadamard derivatives
		and fractional Hyper-Bessel-type operators. Some of the main properties of this class of functions
		and of the related operators will bee analyzed and discussed.
		Then we will show the advantage of this approach to solve some
		integro-differential equations and present a concrete example treating
		a modified Lamb--Bateman equation.

	\section{Basics concepts on Hadamard fractional calculus and Laguerre derivatives}
	
		In order to make the paper self-contained, here we recall some definitions
		and properties of Hadamard fractional integration and differentiation.
		We refer to \citet{Kilbas} for a more detailed analysis. 
      	The first step is to define the Hadamard fractional integral and derivative operators.
		\begin{Definition}
			\label{Integral}
			Let $\Re (\alpha) >0$. The Hadamard fractional integral of order $\alpha$,
			applied to the function $f \in L^p[a,b]$, $0 \le a<b \le \infty$, $x \in [a,b]$,
			is defined as
			\begin{equation}
				J^{\alpha} f(x) =  \frac{1}{\Gamma(\alpha)}
				\int_a^{x} \left( \log \frac{x}{t }\right)^{\alpha-1} f(t) \frac{\mathrm dt}{t}. 
			\end{equation}	
		\end{Definition}
		Before constructing the corresponding derivative operator, we define the following
		space of functions.
		\begin{Definition}
			Let $[a,b]$ be a finite interval such that $-\infty < a<b < \infty$ and let $AC[a,b]$ be
			the space of absolutely continuous functions on $[a,b]$. Let us denote $\delta= x\frac{\mathrm d}{\mathrm dx}$
			and define the space
			\begin{align}
				AC_\delta^n[a,b] = \left\{ g \colon [a,b] \rightarrow \mathbb{C}
				\colon \delta^{n-1} \left[ g(x) \right] \in AC[a,b] \right\}.
			\end{align}
		\end{Definition}
		Clearly $AC_\delta^1[a,b] \equiv AC[a,b]$.
		Analogously to the Riemann--Liouville fractional calculus, the Hadamard fractional derivative
		is defined in terms of the Hadamard fractional integral in the following way.
		\begin{Definition}
			\label{Derivative}
			Let $\delta = x\frac{\mathrm d}{\mathrm dx}$, $\Re(\alpha)>0$
			and $n = [\alpha]+1$, where $[\alpha]$ is the integer part of $\alpha$. 
			The Hadamard fractional derivative of order $\alpha$ applied to the function
			$f \in AC_\delta^n[a,b]$, $0 \le a<b<\infty$, is defined as
			\begin{equation}
				D^{\alpha}f(x)= \delta^n \left(J^{n-\alpha}f \right) (x).
			\end{equation}	
		\end{Definition} 
		A more general Hadamard-type fractional operator has been studied by
		\citet{Kilbas}. We will stick anyway to the classical definition.
		
		It has been proved (see e.g.\ \citet[Theorem 4.8]{Kilbas})
		that in the space $L^p[a,b]$, $0<a<b<\infty$, $1 \le p \le \infty$, the Hadamard fractional
		derivative is the left-inverse operator to the
		Hadamard fractional integral, i.e.\
		\begin{equation}
			\label{comp}
			D^{\alpha}J^{\alpha}f(x)= f(x).
		\end{equation}
		Furthermore, we recall some recent results of \citet{Babusci1} on the operator
		$\left( x\frac{d}{dx} \right)^{-\alpha}$, with 
		$\alpha >0$. The authors have proved by means of operational techniques that it corresponds to the
		Hadamard fractional integral as it is shown in the
		following.
		\citet{Babusci1} have considered the following integral representation of the power $\alpha>0$
		of a general operator $a$
		\begin{equation}
			a^{-\alpha} = \frac{1}{\Gamma(\alpha)} \int_0^{\infty} \mathrm ds \, e^{-sa}s^{\alpha-1},
		\end{equation}	
		so that
		\begin{equation}
			\label{ba}
			\left( x \frac{\mathrm d}{\mathrm dx} \right)^{-\alpha}
			=\frac{1}{\Gamma(\alpha)} \int_0^{\infty} \mathrm ds \, e^{-sx \frac{\mathrm d}{\mathrm dx}}
			s^{\alpha-1}.
		\end{equation}
		To explain it clearly we recall the properties of the dilation operator
		$e^{\lambda x \frac{\mathrm d}{\mathrm dx}}$ \citep{Ricci1}. 
		By considering $x= e^{\theta}$ we have
		\begin{align}
			e^{\lambda x\frac{\mathrm d}{\mathrm dx}} f(x)
			= e^{\lambda \frac{\mathrm d}{\mathrm d \theta}}f (e^{\theta} )
			= f(e^{\theta+\lambda})=f(e^{\lambda}x),
		\end{align}
		where we used the shift action of the exponential operator, that is
		$e^{\lambda \frac{\mathrm d}{\mathrm d\theta}} f(\theta) = f(\theta+\lambda)$.
		Going back to \eqref{ba}
		\begin{align}
			\left( x\frac{\mathrm d}{\mathrm dx} \right)^{-\alpha}
			f(x) & = \frac{1}{\Gamma(\alpha)} \int_0^{\infty} \mathrm ds \, e^{-sx
			\frac{\mathrm d}{\mathrm dx}} f(x) s^{\alpha-1} \\
			& = \frac{1}{\Gamma(\alpha)} \int_0^{\infty} \mathrm ds \, f(e^{-s}x) s^{\alpha-1} \notag \\
			& = \frac{1}{\Gamma(\alpha)} \int_0^{x} \left( \log \frac{x}{t}\right)^{\alpha-1}
			f(t) \frac{\mathrm dt}{t} = J^{\alpha} f(x).
			\notag
		\end{align}

		Since \eqref{comp} holds we immediately have that $\left( x\frac{\mathrm d}{\mathrm dx} \right)^\alpha
		\equiv D^\alpha$.
		This is consistent with the interpretation of the Hadamard fractional derivative
		as a fractional power of the operator $\delta$ \citep[Section 18.3]{libro}.

		These considerations, that  will be clear in the following , imply 
		for example that the Lamb--Bateman equation (see Section \ref{lamb}) can be restated
		as an equation involving a Hadamard fractional derivative.
	
		We now recall the definition of the Laguerre derivative (see e.g.\ \citet{Dattoli,Dattoli1}) as the operator
		$ D_L = \frac{\mathrm d}{\mathrm dx} x \frac{\mathrm d}{\mathrm dx} = \frac{\mathrm d}{\mathrm dx} \delta$.
		Furthermore, let us denote with $D_{nL} = \frac{\mathrm d}{\mathrm dx} x \frac{\mathrm d}{\mathrm dx}
		\dots \frac{\mathrm d}{\mathrm dx} x \frac{\mathrm d}{\mathrm dx}$, $n \in \mathbb{N} \cup \{ 0 \}$,
		the Laguerre operator of order $n$ and containing $n+1$ derivatives.
		For example $D_{0L} \equiv \frac{\mathrm d}{\mathrm dx}$,
		$D_{1L} \equiv D_L$, and $D_{2L} = \frac{\mathrm d}{\mathrm dx} x \frac{\mathrm d}{\mathrm dx}
		x \frac{\mathrm d}{\mathrm dx}$. 
		For the subsequent developments we also recall the following theorem
		by \citet[Theorem 2.2]{Dattoli1} concerning the eigenfunctions of the Laguerre operator
		of order $n$.
		\begin{Theorem}
			The L-exponential function of order $n$
			\begin{equation}
				\label{lexp}
				e_n(x)=\sum_{k = 0}^{\infty}\frac{x^k}{k!^{n+1}}, \qquad n \in \mathbb{N} \cup \{ 0 \},
				\: x \in \mathbb{R},
			\end{equation}
			is an eigenfunction of the operator $D_{nL}$, i.e.\ $D_{nL}e_n(ax)= a\, e_n(ax)$, $a \in \mathbb{R}$.
		\end{Theorem}
		The L-exponential function is also called in literature as $nL$-exponentials
		\citep{Ricci1} and includes the two relevant specific cases of the classical exponential function
		($n=0$) and the $0$th-order Bessel--Tricomi function ($n=1$).
 		\begin{Definition}
 			\label{Derivative1}
			Let $\alpha\in(0,\infty)$, $x > 0$, we define the following
			operator of order $\alpha$.
			\begin{equation}
				\label{aaope}
				\mathfrak{D}^\alpha f(x) = \frac{\mathrm d}{\mathrm dx}D^{\alpha}f(x).
			\end{equation}	
		\end{Definition} 
		\begin{Remark}
		 	Notice that when $\alpha \in (0,1)$, the operator \eqref{aaope} can be written
		 	as
		 	\begin{align}
		 		\mathfrak{D}^\alpha f(x) = \frac{\mathrm d}{\mathrm dx}
				\left(x \frac{\mathrm d}{\mathrm dx} \right) \left(J^{1-\alpha}f \right) (x)		 		
		 		= D_L \left(J^{1-\alpha}f \right)(x),
		 	\end{align}
		 	and thus it can be seen as a generalization of Hadamard derivatives by means of Laguerre
		 	derivatives.
		\end{Remark}

		In the following section we will show some advantages of this approach in order to solve explicitly
		a class of integro-differential equations involving Hadamard and Laguerre type derivatives.
	
	\section{Equations involving the operator $\bm{\mathfrak{D}^\alpha}$}
	
		We first introduce a rather general function which will be specialized in several ways.
		\begin{Definition}[$\alpha$-Mittag--Leffler function]
			The $\alpha$-Mittag--Leffler function $E_{\alpha;\nu,\gamma}(x)$, $x \in \mathbb{R}$, $\alpha>-1$,
			$\nu>0$, $\gamma \in \mathbb{R}$, is defined as
			\begin{align}
				\label{marrMitt}
				E_{\alpha;\nu,\gamma}(x) = \sum_{k=0}^\infty \frac{x^k}{\Gamma^{\alpha+1}(\nu k + \gamma)}.
			\end{align}
		\end{Definition}
		Clearly,
		for $\alpha=0$ the classical Mittag--Leffler function
		$E_{0;\nu,\gamma}(x) = E_{\nu,\gamma}(x)$ is retrieved.
		The following theorem presents the associated Laplace transform for $\alpha \in (0,\infty)$,
		$x \ge 0$.
		\begin{Theorem}
			We have that
			\begin{align}
				\label{mitta-l}
				\int_0^\infty e^{-sx} x^{\gamma-1} E_{\alpha;\nu,\gamma}(\lambda x^\nu) \, \mathrm dx
				= \frac{1}{s^\gamma} E_{\alpha-1;\nu,\gamma}\left( \frac{\lambda}{s^\nu} \right),
				\qquad s > 0, \: \lambda \in \mathbb{R}, \: x \ge 0, \: \alpha \in (0,\infty).
			\end{align}
			\begin{proof}
				By direct calculation we obtain
				\begin{align}
					& \int_0^\infty e^{-sx} x^{\gamma-1} E_{\alpha;\nu,\gamma}(\lambda x^\nu) \, \mathrm dx
					= \sum_{k=0}^\infty \frac{\lambda^k}{\Gamma^{\alpha+1}(\nu k + \gamma)}
					\int_0^\infty e^{-s x} x^{\nu k +\gamma -1} \mathrm dx \\
					& = \sum_{k=0}^\infty \frac{\lambda^k}{\Gamma^{\alpha+1}(\nu k + \gamma)}
					\frac{\Gamma(\nu k + \gamma)}{s^{\nu k+\gamma}}
					= \frac{1}{s^\gamma} E_{\alpha-1;\nu,\gamma}\left( \frac{\lambda}{s^\nu} \right),
					\notag
				\end{align}
				where the inversion of the sum with the integral is permitted by Theorem 30.1 of \citet{doetsch}.
			\end{proof}
		\end{Theorem}
		Note that when $\alpha \rightarrow 0$ we obtain, as expected, the Laplace transform of the
		two-parameters Mittag--Leffler function, i.e.\
		\begin{align}
			\int_0^\infty e^{-sx} x^{\gamma-1} E_{\nu,\gamma}(\lambda x^\nu) \, \mathrm dx
			= \frac{s^{\nu-\gamma}}{s^{\nu}-\lambda}, \qquad s > |\lambda|^{1/\nu}.
		\end{align}
		
		\begin{Theorem}
			For $\alpha=n \in \mathbb{N} \cup \{ 0 \}$, the $n$-Mittag--Leffler function $E_{n;\nu,1}
			(\lambda x^\nu)$, $\lambda \in \mathbb{R}$, $x \ge 0$, $\nu > 0$, is an eigenfunction
			of a fractional hyper-Bessel-type operator
			\begin{align}
				\label{caputoite}
				\underbrace{\frac{\mathrm d^\nu}{\mathrm dx^\nu} x^\nu \frac{\mathrm d^\nu}{\mathrm dx^\nu}
				\dots \frac{\mathrm d^\nu}{\mathrm dx^\nu} x^\nu
				\frac{\mathrm d^\nu}{\mathrm dx^\nu}}_{\text{$n+1$ derivatives}},
			\end{align}
			where $d^\nu/dx^\nu$ represents in this case the Caputo fractional derivative \citep{podlubny}.
			\begin{proof}
				The statement is readily proved by simply observing that
				\begin{align}
					& \underbrace{\frac{\mathrm d^\nu}{\mathrm dx^\nu} x^\nu \frac{\mathrm d^\nu}{\mathrm dx^\nu}
					\dots \frac{\mathrm d^\nu}{\mathrm dx^\nu} x^\nu
					\frac{\mathrm d^\nu}{\mathrm dx^\nu}}_{\text{$n+1$ derivatives}}
					\sum_{k=0}^\infty \frac{\lambda^k x^{\nu k}}{\Gamma^{n+1}(\nu k + 1)} \\
					& = \underbrace{\frac{\mathrm d^\nu}{\mathrm dx^\nu} x^\nu \frac{\mathrm d^\nu}{\mathrm dx^\nu}
					\dots \frac{\mathrm d^\nu}{\mathrm dx^\nu} x^\nu
					\frac{\mathrm d^\nu}{\mathrm dx^\nu}}_{\text{$n$ derivatives}}
					\sum_{k=1}^\infty \frac{\lambda^k x^{\nu k}}{\Gamma^{n+1}(\nu k + 1)}
					\frac{\Gamma(\nu k+1)}{\Gamma(\nu k+1-\nu)}
					\notag \\
					& = \dots = \frac{\mathrm d^\nu}{\mathrm dx^\nu}
					\sum_{k=1}^\infty \frac{\lambda^k x^{\nu k}}{\Gamma^{n+1}(\nu k + 1)}
					\frac{\Gamma^n(\nu k+1)}{\Gamma^n(\nu k+1-\nu)} \notag \\
					& = \sum_{k=1}^\infty
					\frac{\lambda^k x^{\nu k -\nu}}{\Gamma^{n+1}(\nu k +1-\nu)}
					= \lambda E_{n;\nu,1}(\lambda x^\nu), \qquad \lambda \in \mathbb{R}, \: x \ge 0. \notag
				\end{align}
			\end{proof}
		\end{Theorem}
		
		\begin{Remark}
			Notice that if $\alpha = n \in \mathbb{N} \cup \{ 0 \}$ we can write the $\alpha$-Mittag--Leffler
			function \eqref{marrMitt} as a generalized Wright function. Indeed
			\begin{align}
				E_{n;\nu,\gamma}(x) = \sum_{k=0}^\infty \frac{x^k}{\Gamma^{n+1}(\nu k + \gamma)}
				= {}_1 \psi_{n+1}
				\left[ x \left|
				\begin{array}{l}
					(1,1) \\
					\underbrace{(\gamma,\nu),\dots,(\gamma,\nu)}_{\text{$n+1$ times}}
				\end{array}
				\right. \right], \qquad x \in \mathbb{R}, \: \nu >0, \: \gamma \in \mathbb{R},
			\end{align}
			where
			\begin{align}
				{}_p \psi_q
				\left[ z \left|
				\begin{array}{l}
					(a_1,A_1),\dots,(a_p,A_p) \\
					(b_1,B_1),\dots,(b_q,B_q)
				\end{array}
				\right. \right]
				= \sum_{k=0}^\infty \frac{\prod_{j=1}^p \Gamma(a_j+A_j k)}{\prod_{j=1}^q \Gamma(b_j+B_j k)}
				\frac{z^k}{k!},
			\end{align}
			with $a_j,b_j \in \mathbb{C}$, $A_j,B_j \in (0,\infty)$, $z \in \mathbb{C}$,
			is the generalized Wright function.
			
			Clearly, the $n$-Mittag--Leffler function admits other representations, for example as a multi-index
			Mittag--Leffler function \citep{kir}.
		\end{Remark}

		\begin{Remark}
			It is worth noticing that for $(\nu,\gamma)=(1,1)$  the operator \eqref{caputoite} and the function
			\eqref{marrMitt} coincide respectively with the operator $D_{nL}$ and the L-exponential function
			\eqref{lexp}.
		\end{Remark}

		By means of the following definition we specialize function \eqref{marrMitt} for $(\nu,\gamma)=(1,1)$
		thus obtaining a generalization of the L-exponential function
		\eqref{lexp} for fractional values of parameter $n$.
		\begin{Definition}
			The $\alpha L$-exponential function $\mathfrak{e}_\alpha (x)$, $\alpha \in (-1,\infty)$, $x \in \mathbb{R}$, is defined
			as
			\begin{equation}
				\label{lhexp}
				\mathfrak{e}_\alpha(x)=\sum_{k = 0}^{\infty}\frac{x^k}{k!^{\alpha+1}}.
			\end{equation}
		\end{Definition}
		For the Cauchy--Hadamard theorem, the $\alpha L$-exponential function is clearly convergent in $\mathbb{R}$ for
		$\alpha \in (-1,\infty)$. Its Laplace transform for $x\ge 0$, $\alpha \in (0,\infty)$ reads
		\begin{align}
			\label{goodwine}
			\int_0^\infty e^{-s x} \mathfrak{e}_\alpha (x)
			\, \mathrm dx = \frac{1}{s} \sum_{k=0}^\infty \frac{s^{-k}}{k!^\alpha}
			= \frac{1}{s} \; \mathfrak{e}_{\alpha-1} \left( \frac{1}{s} \right), \qquad s>0,
		\end{align}
		where the inversion of the sum with the integral is permitted by Theorem 30.1 of \citet{doetsch}.
		Note that formula \eqref{goodwine} for $\alpha \rightarrow 0$, $s > 1$,
		correctly gives the Laplace transform of the
		classical exponential function. Note that, as expected, the above result is consistent with \eqref{mitta-l}.
		\begin{figure}
			\centering
			\includegraphics[scale=.38]{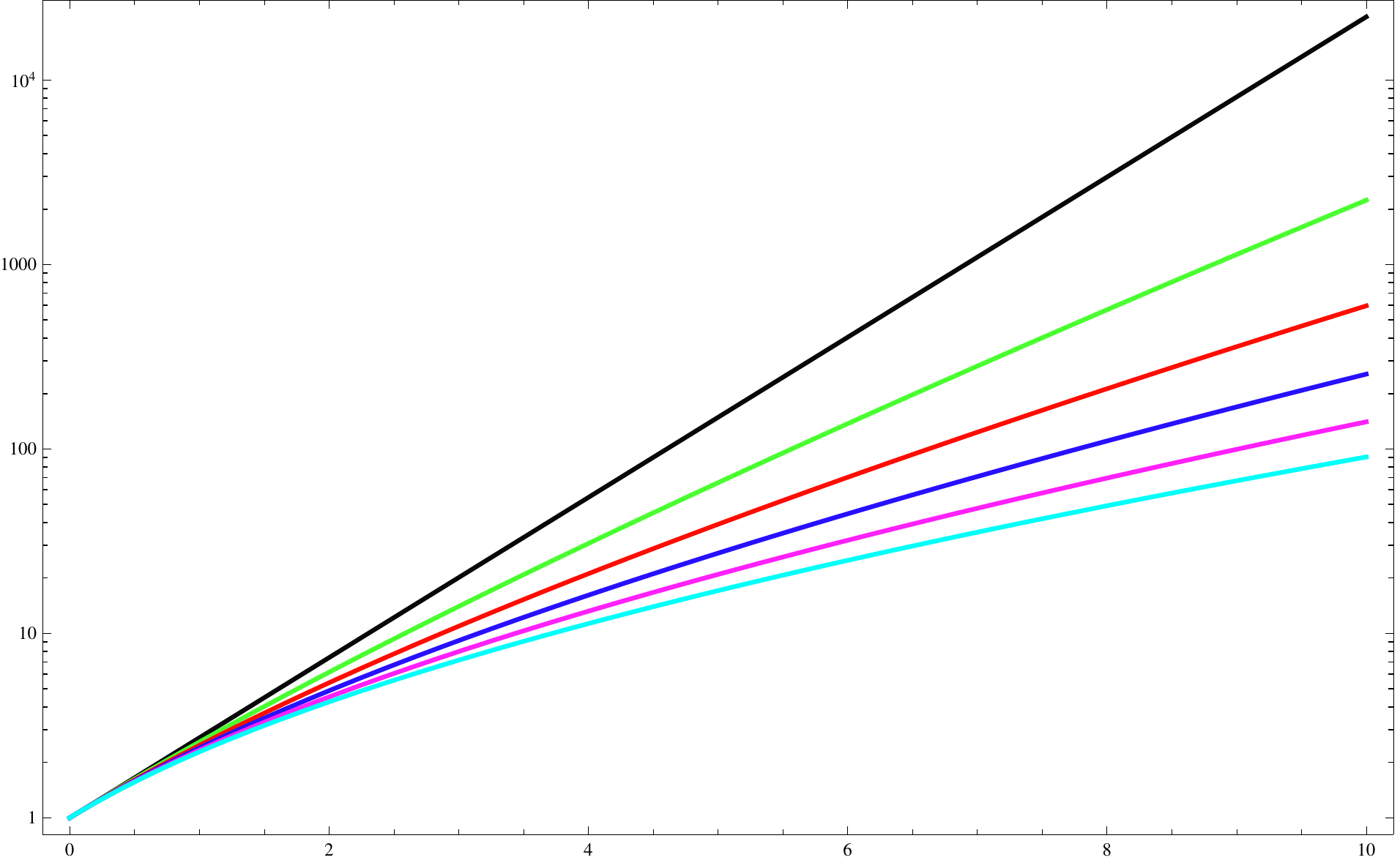}
			\includegraphics[scale=.38]{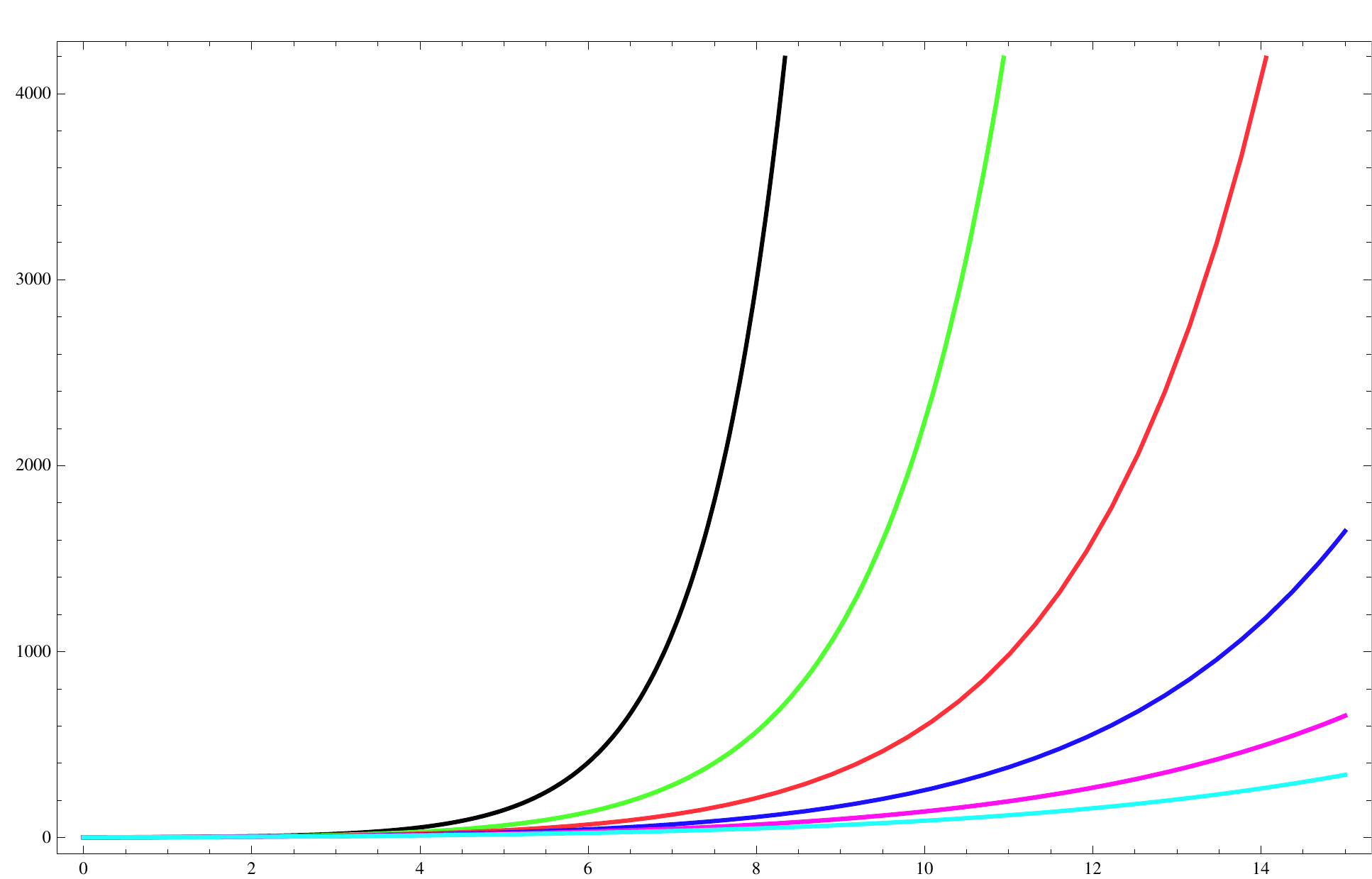}
			\caption{\label{fig}The $\alpha L$-exponential function $\mathfrak{e}_\alpha (x)$,
			$x \ge 0$ in loglinear plot (left)
			and linear plot (right), for $\alpha=0$
			(black---classical exponential), $\alpha=0.2$ (green), $\alpha=0.4$ (red), $\alpha=0.6$ (blue),
			$\alpha=0.8$ (magenta), $\alpha=1$ (cyan---$0$th-order Bessel--Tricomi function).}
		\end{figure}

\begin{figure}
			\centering
			\includegraphics[scale=.6]{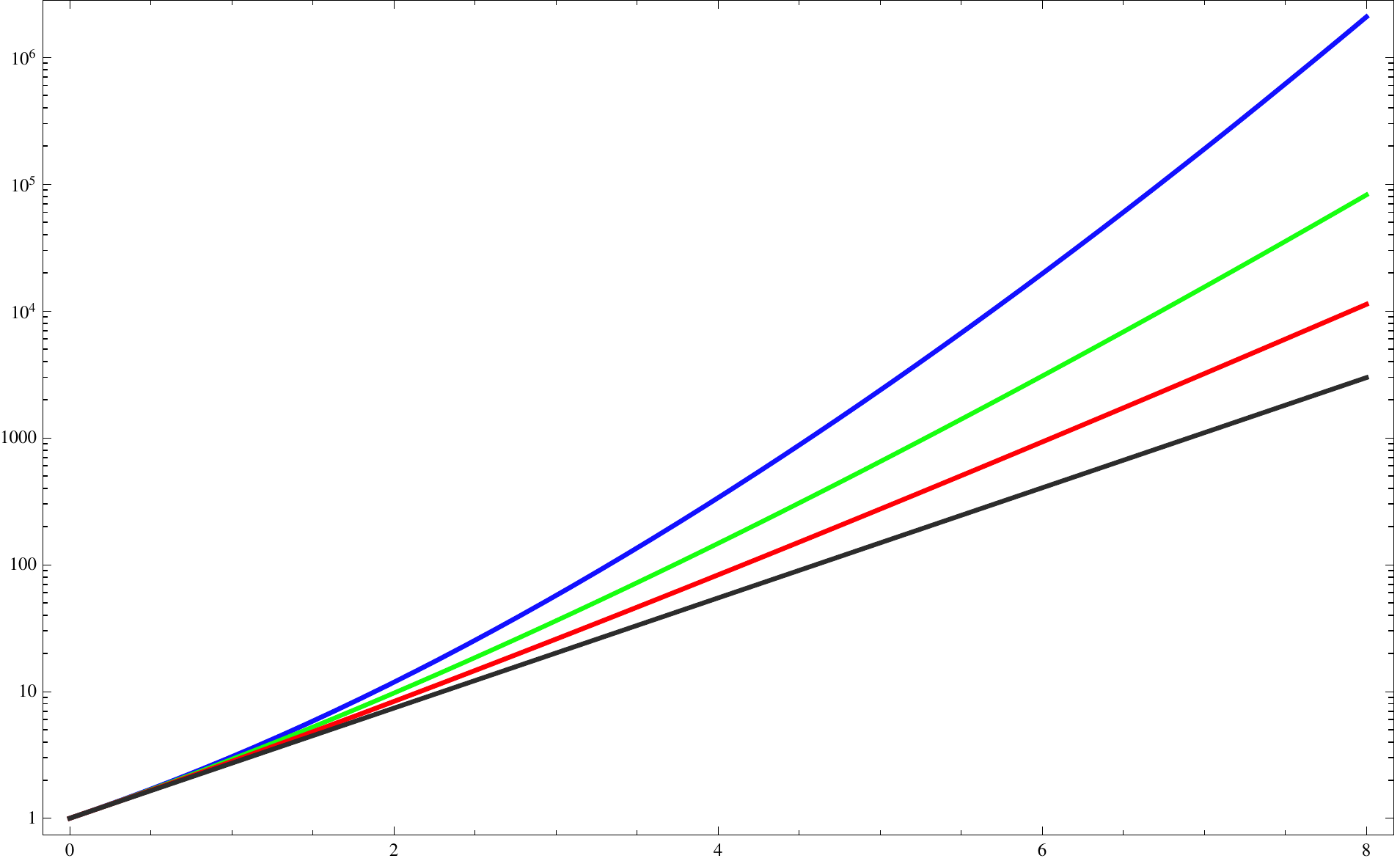}
			\caption{\label{fig2}The $\alpha L$-exponential function $\mathfrak{e}_\alpha(x)$, $x \ge 0$ in loglinear plot
			for $\alpha=0$
			(black---classical exponential), $\alpha=-0.1$ (red), $\alpha=-0.2$ (green), $\alpha=-0.3$ (blue).}
		\end{figure}

		Figures \ref{fig} and \ref{fig2} show the $\alpha L$-exponential function for some values of parameter $\alpha$.

		\begin{Theorem}
			\label{teoremolo}
			The $\alpha L$-exponential function $\mathfrak{e}_\alpha(\lambda x)$, $\alpha \in [0,\infty)$, $x \ge 0$,
			$\lambda \in \mathbb{R}$, is an eigenfunction of the operators
			\begin{align}
				\underbrace{\frac{\mathrm d}{\mathrm dx} x \frac{\mathrm d}{\mathrm dx}
				\dots \frac{\mathrm d}{\mathrm dx} x \frac{\mathrm d}{\mathrm dx}}_{\text{$r$ derivatives}}
				D^{\alpha+1 -r}, \qquad \alpha \in [0,\infty), \: r=1, \dots, n,
			\end{align}
			where $n=[\alpha+1]$ is the integer part of $\alpha+1$. Moreover, for $\alpha \in (-1,0)$
			the $\alpha L$-exponential function $\mathfrak{e}_\alpha(\lambda x)$,
			$x \ge 0$, $\lambda \in \mathbb{R}$, is an
			eigenfunction of the operator $x^{-1} D^{\alpha+1}$.
			
			\begin{proof}
				For the first part we can write that
				\begin{align}
					& \underbrace{\frac{\mathrm d}{\mathrm dx} x \frac{\mathrm d}{\mathrm dx}
					\dots \frac{\mathrm d}{\mathrm dx} x \frac{\mathrm d}{\mathrm dx}}_{\text{$r$ derivatives}}
					D^{\alpha+1 -r} \mathfrak{e}_\alpha(\lambda x)
					= \underbrace{\frac{\mathrm d}{\mathrm dx} x \frac{\mathrm d}{\mathrm dx}
					\dots \frac{\mathrm d}{\mathrm dx} x \frac{\mathrm d}{\mathrm dx}}_{\text{$r$ derivatives}}
					\sum_{k=1}^\infty \lambda^k \frac{k^{\alpha+1-r}x^k}{k!^{\alpha+1}} \\
					& = \frac{\mathrm d}{\mathrm dx} \sum_{k=1}^\infty \lambda^k \frac{k^{\alpha}x^k}{k!^{\alpha+1}}
					= \lambda \mathfrak{e}_\alpha(\lambda x). \notag
				\end{align}
				Analogously, in order to prove the second part it is sufficient to note that
				\begin{align}
					x^{-1} D^{\alpha+1} \mathfrak{e}_\alpha(\lambda x) = x^{-1} \sum_{k=1}^\infty
					\lambda^k \frac{k^{\alpha+1} x^k}{k!^{\alpha+1}} = \lambda \mathfrak{e}_\alpha(\lambda x).
				\end{align}		
			\end{proof}
		\end{Theorem}
		\begin{Remark}
			An interesting specific case of Theorem \ref{teoremolo} is when $r=1$. In this case the
			$\alpha L$-exponential function is an eigenfunction of the operator $\frac{\mathrm d}{\mathrm dx}D^{\alpha}
			= \mathfrak{D}^\alpha$ (see Definition \ref{Derivative1}).
			This can be explained also by recalling that
			\begin{align}
				\underbrace{\frac{\mathrm d}{\mathrm dx} x \frac{\mathrm d}{\mathrm dx}
				\dots \frac{\mathrm d}{\mathrm dx} x \frac{\mathrm d}{\mathrm dx}}_{\text{$r$ derivatives}}
				D^{\alpha+1 -r} = \frac{\mathrm d}{\mathrm dx} \delta^{r-1} D^{\alpha+1-r}
				= \frac{\mathrm d}{\mathrm dx} \delta^{r-1} \delta^{n-r+1} J^{n-\alpha}
				= \frac{\mathrm d}{\mathrm dx} D^\alpha = \mathfrak{D}^\alpha.
			\end{align}
		\end{Remark}

		In order to describe some applications
		we first recall the following result \cite[page 490, Theorem 5.1]{Dattoli1} specialized for $n=1$.
		\begin{Theorem}
			Consider the problem
			\begin{align}
				\begin{cases}
					D_{L,x} \; S(x,t) = \frac{\partial}{\partial t} S(x,t), & x > 0, t \ge 0, \\
					S(0,t) = s(t),
				\end{cases}
			\end{align}
			with the analytic boundary condition $s(t)$. The operational solution is given by
			\begin{equation}
				\label{30ff}
				S(x,t)= e_1 \left(x\frac{\partial}{\partial t}\right) s(t)
				= \sum_{k=0}^{\infty} \frac{x^k}{k!^{2}}\frac{\partial^k}{\partial t^k}s(t).
			\end{equation}
		\end{Theorem}
		The operational solution \eqref{30ff} becomes an effective solution when the series converges,
		and this depends upon the actual form of the boundary condition $S(0,t)$.
		Notice that, from a physical point of view, this \emph{Laguerre-heat equation}
		is in fact a diffusion equation with a space-dependent diffusion coefficient.
		
		In a similar way we can solve more general fractional partial differential equations.
		We notice that the utility of operational methods to solve fractional differential equations was
		pointed out by different authors, see for example \citet{uno,due,tre,quattro}. See also
		the book by \citet{librot} and the references therein.
		
		Here we present the following result for an initial value problem involving the operator $\mathfrak{D}^\alpha$.
		\begin{Theorem}
			\label{30bb}
			Consider the following initial value problem.
			\begin{equation}
				\label{IVP}
				\begin{cases}
					 \mathfrak{D}_t^{\alpha} f(x,t) = \Theta_x f(x,t), & t > 0, \\ 
					 f(x,0) = g(x),
				\end{cases}
			\end{equation}
			for $\alpha \in (0,1)$, with an analytic initial condition $g(x)$
			and where $\Theta_x$ is a generic linear integro-differential operator with constant coefficients,
			acting on $x$, and
			which satisfies the semigroup property, that is $\Theta_x \Theta_y = \Theta_{x+y}$.
			By applying Theorem \ref{teoremolo} for $r=1$, $\alpha \in (0,1)$,
			the operational solution to \eqref{IVP} is given by
			\begin{equation}
				\label{A}
				f(x,t) = \mathfrak{e}_\alpha ( t \, \Theta_x ) \,
				g(x) = \sum_{r=0}^{\infty} \frac{t^r \Theta^r_x}{r!^{\alpha+1}} g(x).
			\end{equation}
		\end{Theorem}
		Again, when the operational solution \eqref{A} converges, it becomes the effective solution.
		\begin{Corollary}
			Consider the following boundary value problem.  
			\begin{equation}
				\label{BVP}
				\begin{cases}
					\mathfrak{D}_x^{\alpha} w(x,t) = \Xi_t w(x,t), & \alpha \in (0,1), \: t \ge 0, \: x > 0, \\
					w(0,t) = h(t),
				\end{cases}
			\end{equation}
			with an analytic boundary condition $h(t)$
			and where $\Xi_t$ is a generic linear integro-differential operator with constant coefficients,
			acting on $t$, and which satisfies the semigroup property. Analogously to Theorem \ref{30bb},
			the operational solution to \eqref{BVP} is given by
			\begin{equation}
				\label{B}
				w(x,t) = \mathfrak{e}_\alpha ( x \, \Xi_t ) \,
				h(t) = \sum_{r=0}^{\infty} \frac{x^r\Xi^r_t}{r!^{\alpha+1}} h(t).
			\end{equation}
		\end{Corollary}
		It is important to remark that the operational solutions
		\eqref{A} and \eqref{B} are indeed effective solutions whenever the corresponding
		series are convergent.

		\begin{ex}
			Consider the following initial value problem.
			\begin{equation}
				\begin{cases}
					 \mathfrak{D}_t^{\alpha} \; u(x,t) = \frac{\partial^2}{\partial x^2} u(x,t),
					 & \alpha \in (0,1), \: t > 0, \\
					 u(x,0) = \sin x, 
				\end{cases}
			\end{equation}
			By Theorem \ref{30bb} the solution is given by
			\begin{equation}
				u(x,t) = \sin x \; \mathfrak{e}_\alpha (t) = \sin x \, \sum_{k = 0}^{\infty}
				\frac{t^k}{k!^{\alpha+1}},
			\end{equation}
			which is clearly convergent for each $(x,t) \in \mathbb{R} \times \mathbb{R}^+$.
		\end{ex}

	\section{An application to a modified Lamb--Bateman equation}
	
		\label{lamb}		
		In \cite{Babusci1} the authors have shown that a modified Lamb--Bateman integral equation
		can be written as a fractional equation involving 
		Hadamard derivatives of order $1/2$. They have studied the equation
		\begin{equation}
			\label{Lamb}
			\int_0^{\infty} \mathrm dy \; u \left( e^{-y^2} x \right)  = f(x),
		\end{equation}
		where $u(x)$ is the unknown function that has to be determined and $f(x)$ is a known continuous function.
		Recalling that, for the dilation operator $\hat{E}(\lambda) = e^{\lambda x \frac{\mathrm d}{\mathrm dx}}$, we have
		\begin{align}
			\hat{E}(\lambda) \, u(x) = e^{\lambda x \frac{\mathrm d}{\mathrm dx}} u(x) = u(e^{\lambda}x),
		\end{align}
		then equation \eqref{Lamb} can consequently be written as 
		\begin{equation}
			\int_0^{\infty} \mathrm dy \, e^{-y^2 x \frac{\mathrm d}{\mathrm dx}} u(x) = f(x),
		\end{equation}
		and therefore
		\begin{equation}
			u(x) = \frac{2}{\sqrt{\pi}} \left( x\frac{\mathrm d}{\mathrm dx} \right)^{1/2} \hspace{-.3cm} f(x)
			= \frac{2}{\sqrt{\pi}} D^{1/2} f(x),
		\end{equation} 
		that is, an integral equation involving Hadamard derivatives of order $1/2$.
		
		Note that, in the general case, the function $f$ depends on $u$. The simplest
		choice is $f(u(x)) = u(x)$. In this specific case we have
		\begin{equation}
			u(x) = \frac{2}{\sqrt{\pi}} D^{1/2} u(x).
		\end{equation}
		An application of the $x$-derivative to both sides of the above equation leads to a formulation of the
		Lamb--Bateman equation involving the operator $\mathfrak{D}^\alpha$.
		\begin{equation}
			\mathfrak{D}^{1/2} u(x) = \frac{\sqrt{\pi}}{2}\frac{\mathrm d}{\mathrm dx} u(x).
		\end{equation}
		By using a simple \textit{ansatz}, that is $u(x)= x^\beta$, where $\beta$ is to be specified,
		we find by substitution
		\begin{equation}
			\beta^{3/2}x^{\beta-1}=\frac{\sqrt{\pi}}{2}\beta x^{\beta-1},
		\end{equation}
		so that $u(x)= x^{\beta}$ is a simple exact solution if $\beta=\pi/4$.
		
		In general we can consider the equation
		\begin{align}
			\int_0^\infty \mathrm dy \, u \left( e^{-y^{1/\mu}} x \right) = f(x), \qquad \mu \in (0,\infty),
		\end{align}
		which can be rewritten, as before, as
		\begin{align}
			\int_0^\infty \mathrm dy \, e^{-y^{1/\mu}x\frac{\mathrm d}{\mathrm dx}} u(x) = f(x).
		\end{align}
		By evaluating the integral we then have that
		\begin{align}
			u(x) = \frac{1}{\Gamma(\mu+1)} \left( x\frac{\mathrm d}{\mathrm dx} \right)^\mu f(x)
			= \frac{1}{\Gamma(\mu+1)} D^\mu f(x).
		\end{align}
		Again, the simplest choice of $f(u(x)) = u(x)$ leads to
		\begin{align}
			u(x) = \frac{1}{\Gamma(\mu+1)} D^\mu u(x), \qquad \mu \in (0,\infty).
		\end{align}
		The related equation involving the operator $\mathfrak{D}^\alpha$ is
		\begin{align}
			\mathfrak{D}^\mu u(x) = \Gamma(\mu+1) \frac{\mathrm d}{\mathrm dx} u(x),
		\end{align}
		and by choosing again $u(x)=x^\beta$ we obtain that $\beta = \left[ \Gamma(\mu+1) \right]^{1/\mu}$.

\end{document}